\documentclass[12pt,a4paper]{article}

\usepackage[english]{babel}
\usepackage{amssymb}
\usepackage{amsmath}
\usepackage[pdftex]{graphicx}
\usepackage{graphics}
\usepackage{times}
\usepackage{url}
\usepackage{multicol}

\newtheorem{tm}{Theorem}
\newtheorem{defi}{Definition}

\newcommand{\linesp}[1]{%
\baselineskip #1\normalbaselineskip plus 0.03\normalbaselineskip 
                minus 0.02\normalbaselineskip 
\parskip .35\baselineskip plus 0.03\normalbaselineskip 
                minus 0.03\normalbaselineskip 
}

\begin{document}

 \begin{center} {\large {\bf Generalized Rose Surfaces and their Visualizations }}

\bigskip

 SONJA GORJANC\\
{\small University of Zagreb, Faculty of Civil Engineering,\\ Ka\v{c}i\'{c}eva 26, 10000 Zagreb, Croatia\\
e-mail: sgorjanc@master.grad.hr}

\medskip

 EMA JURKIN\\
{\small University of Zagreb, Faculty of Mining, Geology and Petroleum Engineering,\\ Pierottijeva 6, 10000 Zagreb, Croatia\\
e-mail: ema.jurkin@rgn.hr}
\end{center}

\bigskip

\noindent {\small {\bf Abstract}.\\
In this paper we construct a new class of algebraic surfaces in three-dimensional Euclidean space generated by a cyclic-harmonic curve  and a congruence of circles. We study their properties and visualize them with the program {\it Mathematica}.

\medskip
\medskip

}\\

\section{Introduction}

In  article \cite{Izumiya}  the definition of a {\it circular surface} is given as follows:
\begin{defi}
A {\it circular surface} is (the image of) a map $V:I\times \mathbb R/2\pi \mathbb Z \longrightarrow \mathbb R^3$ defined by
\begin{equation}\label{cs}
V(t,\theta)=\gamma(t)+r(t)(\cos \theta\mathbf a_1(t)+\sin \theta\mathbf a_2(t)),\nonumber
\end{equation}
where $\gamma, \mathbf a_1, \mathbf a_2 \colon I\longrightarrow \mathbb R^3$ and $r \colon I\longrightarrow \mathbb R_{>0}$.
\end{defi}
It is assumed that $\langle \mathbf a_1,\mathbf a_1\rangle =\langle \mathbf a_2,\mathbf a_2\rangle=1$ and $\langle \mathbf a_1,\mathbf a_2\rangle=0$ for all $t\in I$, where $\langle ,\rangle$ denotes the canonical inner product on $R^3$. The curve $\gamma$ is called   a {\it base curve} and the pair of curves $\mathbf a_1$, $\mathbf a_2$ is said to be   a {\it director frame}. The standard circles $\theta\mapsto r(t)(\cos\theta\mathbf a_1(t)+\sin\theta\mathbf a_2(t))$ are called {\it generating circles}.   In this paper the authors classify and investigate the differential properties of this one-parameter family of standard circles with a fixed radius.

In paper \cite{circular} the authors defined three types (elliptic, parabolic and hyperbolic) of congruences  $\mathcal C(p)$ that contain circles passing through two points $P_1$ and $P_2$.  Then, for a given congruence $\mathcal C(p)$  of a certain type and a given curve $\alpha$  they defined a circular surface $\mathcal {CS}(\alpha,p)$  as the system of circles from $\mathcal C(p)$  that intersect $\alpha$. These surfaces contain the generating circles of variable radii.

The {\it rose surfaces} studied in \cite{roses} are circular surfaces $\mathcal {CS}(\alpha,p)$ where $\alpha$  is a rose (rhodonea curve) given by the polar equation $\rho=\cos \frac{n}{d}\varphi$, and  $\mathcal C(p)$ is an elliptic or parabolic congruence. The rose lies in the plane perpendicular to the axis $z$ having the directing point $P_i$ as its multiple point.

In this paper we extend $\alpha$ to all cyclic-harmonic curves $\rho=\cos \frac{n}{d}\varphi+a$, $a\in \mathbb R^+\cup \{0\}$ (see \cite{Hilton/clanak}, \cite {Moritz}),  and include a  hyperbolic congruence  $\mathcal C(p)$.  For such circular surfaces $\mathcal {CS}(\alpha,p)$,  we give an overview of their algebraic properties and visualize their numerous forms with the program {\it Mathematica}.

\section{Cyclic-Harmonic Curve}

A {\it cyclic-harmonic curve}  is given by the following polar equation
\begin{equation}\label{polar}
r(\varphi )= b \cos\frac{n}{d}\varphi +a,\,\varphi\in [0,2d\pi],
\end{equation}
where $\frac{n}{d}$ is a positive rational number in lowest terms and $a, b\in\mathbb R^+$.
 It is a locus of a composition of two simultaneous motions: a {\it simple harmonic} motion $r(\varphi)=  b \cos nt+ a$, and an {\it uniform angular} motion $\varphi =dt$.

\begin{figure}[h]
  \centering
 \includegraphics[scale=0.65]{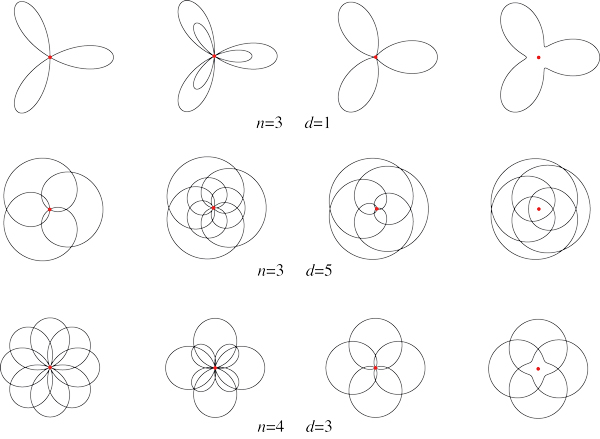}
  \caption{Some examples of $CH(n,d,a)$}
  \label{fig1}
\end{figure}

According to \cite{Moritz}, a cyclic-harmonic curve is called {\it foliate}, {\it prolate}, {\it cuspitate} or {\it curtate} if $a=0$,  $a<b$, $a=b$ or  $a>b$, respectively. Without the loss of generality, we will suppose that $b=1$ and denote a cyclic-harmonic curve by  $CH(n,d,a)$. Some examples of these curves are shown in Figure \ref{fig1}.

\begin{center}


\begin{multicols}{4}

Foliate\\
$a=0$

 Prolate \\
$a<1$

Cuspitate\\
$a=1$  

 Curtate\\
$a>1$
\end{multicols}
\end{center}

The curves $CH(n,d,0)$ are called {\it roses} or {\it rhodonea curves} and have been treated in \cite{roses}. The reader is invited to read the paper and find more about their algebraic properties. Here we give only their implicit equation:
{\small
\begin{equation}\label{rose}
\left( \sum_{k=0}^{\lfloor \frac{d}{2}\rfloor}\sum_{j=0}^{k}(-1)^{j+k} \binom{d}{2k} \binom{k}{j} (x^2+y^2)^{\frac{n+d}{2}-k+j} \right)^s-\left( \sum_{i=0}^{\lfloor \frac{d}{2}\rfloor}(-1)^{i} \binom{n}{2i} x^{n-2i}y^{2i}\right)^s= 0,
\end{equation}
}
where $s=1$ if $n \cdot d$ is odd and $s=2$ if $n \cdot d$ is even.

\begin{tm}\label{rose}

For a cyclic-harmonic curve $CH(n,d,a)$, the following table is valid:

\begin{center}
\rm

\begin{tabular}{|l |c| c |c |c | c| }
\hline

 \multicolumn{2}{|c|}{ } &  \hspace{0.1cm} order of  \hspace{0.1cm}  & \hspace{0.1cm}  multiplicity of \hspace{0.1cm}  &  \hspace{0.1cm}   multiplicity of\\
 \multicolumn{2}{|c|}{$CH(n,d,a)$}& \hspace{0.1cm}  $CH(n,d,a)$  \hspace{0.1cm}   & \hspace{0.1cm}   the origin $O$ \hspace{0.1cm}   & \hspace{0.1cm}   the absolute points  
\\

\hline
\hline
$a=0$ \hspace{0.1cm}  & \hspace{0.1cm}  $d<n$ \hspace{0.1cm}  &  \hspace{0.1cm}  $n+d$  \hspace{0.1cm}  & \hspace{0.1cm}   $n$  \hspace{0.1cm}  &  \hspace{0.1cm}  $d$\\

$n\cdot d$ \hspace{0.1cm}  odd \hspace{0.1cm}  & \hspace{0.1cm}  $d>n$ \hspace{0.1cm}  & \hspace{0.1cm}  $n+d$ \hspace{0.1cm}  & \hspace{0.1cm}  $n$ \hspace{0.1cm}  & \hspace{0.1cm}  $\frac{n+d}{2}$\\
\hline
other \hspace{0.1cm} &  \hspace{0.1cm}  $d<n$ \hspace{0.1cm}  & \hspace{0.1cm} $2(n+d)$ \hspace{0.1cm} & \hspace{0.1cm} $2n$ \hspace{0.1cm} & \hspace{0.1cm} $2d$\\

cases \hspace{0.1cm} & \hspace{0.1cm}  $d>n$ \hspace{0.1cm} &  \hspace{0.1cm} $2(n+d)$ \hspace{0.1cm} & \hspace{0.1cm} $2n$ \hspace{0.1cm} & \hspace{0.1cm} $(n+d)$\\
\hline
\end{tabular}

\bigskip

 Table1: Properties of $CH(n,d,a)$

\end{center}

\end{tm}

\noindent {\bf Proof.} Statements regarding roses were presented in \cite{roses}. So only the cases for $a \neq 0$ are left to be proved.

By substituting the standard formulas for  conversion from the polar to the Cartesian coordinates $
r(\varphi)=\sqrt{x^2+y^2},\,\,\cos\varphi=\frac{x}{\sqrt{x^2+y^2}},\,\,\sin\varphi=\frac{y}{\sqrt{x^2+y^2}}
$
into (\ref{polar}), we obtain the following condition for the points on   $CH(n,d,a)$:
\begin{equation}
\label{cossin}
\cos\frac{n}{d}\varphi=\sqrt{x^2+y^2}-a,\,\,\,\,\sin\frac{n}{d}\varphi=\sqrt{1-(\sqrt{x^2+y^2}-a)^2}.
\end{equation}
Since $
\cos d\frac{n}{d}\varphi=\cos n\varphi,
$
the  cosine multiple angle formula\\
$
\cos n\varphi=\sum_{i=0}^{\lfloor \frac{n}{2}\rfloor} (-1)^i\binom{n}{2i} (\sin\varphi)^{2i}(\cos\varphi)^{n-2i}
$
gives us the following equality:
{\small
\begin{equation}
\label{identity2}
\sum_{j=0}^{\lfloor \frac{d}{2}\rfloor} (-1)^j\binom{d}{2j} (\sin\frac{n}{d}\varphi)^{2j}(\cos\frac{n}{d}\varphi)^{d-2j}=
\sum_{i=0}^{\lfloor \frac{n}{2}\rfloor} (-1)^i\binom{n}{2i} (\sin\varphi)^{2i}(\cos\varphi)^{n-2i}.
\end{equation}
}
Thus,
{\small
\begin{align}
\label{identity3}
&(x^2+y^2)^{\frac{n}{2}}\sum_{j=0}^{\lfloor \frac{d}{2}\rfloor} (-1)^j\binom{d}{2j} (1-(\sqrt{x^2+y^2}-a)^2)^{j}(\sqrt{x^2+y^2}-a)^{d-2j}\nonumber\\
&=\sum_{i=0}^{\lfloor \frac{n}{2}\rfloor} (-1)^i\binom{n}{2i} x^{n-2i}y^{2i}.
\end{align}
}
After applying  the binomial theorem to the expression $(1-(\sqrt{x^2+y^2}-a)^2)^{j}$, equation  (\ref{identity3}) turns into
{\small
\begin{equation}
\label{identity4}
(x^2+y^2)^{\frac{n}{2}}\sum_{j=0}^{\lfloor \frac{d}{2}\rfloor}\sum_{k=0}^{j} (-1)^{2j-k}\binom{d}{2j}\binom{j}{k} (\sqrt{x^2+y^2}-a)^{d-2k}=
\sum_{i=0}^{\lfloor \frac{n}{2}\rfloor} (-1)^i\binom{n}{2i} x^{n-2i}y^{2i}.
\end{equation}
}
Using the same formula again, this time applying it to $ (\sqrt{x^2+y^2}-a)^{d-2k}$, the equation above becomes:
{\small
\begin{align}
\label{identity5}
&(x^2+y^2)^{\frac{n}{2}}\sum_{j=0}^{\lfloor \frac{d}{2}\rfloor}\sum_{k=0}^{j}\sum_{l=0}^{d-2k} (-1)^{d-k-l}\binom{d}{2j}\binom{j}{k}\binom{d-2k}{l} a^{d-2k-l}(x^2+y^2)^{\frac{l}{2}}\nonumber\\
&= \sum_{i=0}^{\lfloor \frac{n}{2}\rfloor} (-1)^i\binom{n}{2i} x^{n-2i}y^{2i}.
\end{align}
}
The left side of equation (\ref{identity5}) can be written in the following form:
{\small
\begin{align}
\label{rastav}
&(x^2+y^2)^{\frac{n}{2}}\sum_{j=0}^{\lfloor \frac{d}{2}\rfloor}\sum_{k=0}^{j}\sum_{l=0}^{\lfloor \frac{d-2k}{2}\rfloor} (-1)^{d-k}\binom{d}{2j}\binom{j}{k}\binom{d-2k}{2l} a^{d-2k-2l}(x^2+y^2)^{l}\nonumber\\
+&(x^2+y^2)^{\frac{n+1}{2}}\sum_{j=0}^{\lfloor \frac{d}{2}\rfloor}\sum_{k=0}^{j}\sum_{l=0}^{\mathrm D} (-1)^{d-k-1}\binom{d}{2j}\binom{j}{k}\binom{d-2k}{2l+1} a^{d-2k-2l-1}(x^2+y^2)^{l},
\end{align}
}
where
\begin{equation}
\label{D}
\mathrm D =
\begin{cases}
\lfloor \frac{d-2k}{2}\rfloor & \text{if {\it d} is an odd number},\\
\lfloor \frac{d-2k}{2}\rfloor -1 & \text{if {\it d} is an even number.}
\end{cases}
\end{equation}
Therefore, the implicit equations of cyclic-harmonic curves are:

\smallskip
\noindent -- if $n$ is an even number ($d$ must be odd)
{\small
\begin{align}
\label{implicit_n_par}
&(x^2+y^2)^{n+1}\Big( \sum_{j=0}^{\lfloor \frac{d}{2}\rfloor}\sum_{k=0}^{j}\sum_{l=0}^{\lfloor \frac{d-2k}{2}\rfloor} (-1)^{d-k-1}\binom{d}{2j}\binom{j}{k}\binom{d-2k}{2l+1} a^{d-2k-2l-1}(x^2+y^2)^{l}\Big)^2\nonumber\\
=&\Big( \sum_{i=0}^{\lfloor \frac{n}{2}\rfloor} (-1)^i\binom{n}{2i} x^{n-2i}y^{2i}\\
&-(x^2+y^2)^{\frac{n}{2}}\sum_{j=0}^{\lfloor \frac{d}{2}\rfloor}\sum_{k=0}^{j}\sum_{l=0}^{\lfloor \frac{d-2k}{2}\rfloor} (-1)^{d-k}\binom{d}{2j}\binom{j}{k}\binom{d-2k}{2l} a^{d-2k-2l}(x^2+y^2)^{l}\Big) ^2\nonumber
\end{align}}
-- if $n$ is an odd number
{\small
\begin{align}
\label{implicit_n_nep}
&(x^2+y^2)^{n}\Big( \sum_{j=0}^{\lfloor \frac{d}{2}\rfloor}\sum_{k=0}^{j}\sum_{l=0}^{\lfloor \frac{d-2k}{2}\rfloor} (-1)^{d-k}\binom{d}{2j}\binom{j}{k}\binom{d-2k}{2l} a^{d-2k-2l}(x^2+y^2)^{l}\Big)^2\nonumber\\
=&\Big( \sum_{i=0}^{\lfloor \frac{n}{2}\rfloor} (-1)^i\binom{n}{2i} x^{n-2i}y^{2i}\\
&- (x^2+y^2)^{\frac{n+1}{2}}\sum_{j=0}^{\lfloor \frac{d}{2}\rfloor}\sum_{k=0}^{j}\sum_{l=0}^{\mathrm D} (-1)^{d-k-1}\binom{d}{2j}\binom{j}{k}\binom{d-2k}{2l+1} a^{d-2k-2l-1}(x^2+y^2)^{l}\Big)^2\nonumber
\end{align}
}
where $\mathrm D$ is given by (\ref{D}).

Since both equations (\ref{implicit_n_par}) and (\ref{implicit_n_nep}) are of the degree $2(n+d)$, the algebraic curve $CH(n,d,a)$ is of the order $2(n+d)$.

 According to \cite[p.251]{Harris}, for an $n$th order  algebraic plane curve which passes through the origin $O(0,0)$ and is given by an algebraic equation $P^n(x,y)=0$, the tangent lines at the origin are given by a homogeneous algebraic equation $f^k(x,y)=0$, where $f^k(x,y)$ is the term of $P^n(x,y)$ with the lowest degree. Thus, $O(0,0)$ is an $2n$-fold point of $CH(n,d,a)$ because the tangent lines are given by the following algebraic equations of degree $2n$:

\bigskip

-- if $n$ is an even number ($d$ must be odd)
\begin{equation}
\label{tangent_n_par}
\Big( \sum_{i=0}^{\lfloor \frac{n}{2}\rfloor} (-1)^i\binom{n}{2i} x^{n-2i}y^{2i}-\mathrm A\cdot (x^2+y^2)^{\frac{n}{2}}\Big) ^2=0
\end{equation}

-- if $n$ is an odd number
{\small
\begin{equation}
\label{tangent_n_nep}
\mathrm A^2\cdot (x^2+y^2)^{n}-\Big( \sum_{i=0}^{\lfloor \frac{n}{2}\rfloor} (-1)^i\binom{n}{2i} x^{n-2i}y^{2i}\Big)^2=0
\end{equation}
}
where
{\small\begin{equation*}
\mathrm A=\sum_{j=0}^{\lfloor \frac{d}{2}\rfloor}\sum_{k=0}^{j} (-1)^{d-k}\binom{d}{2j}\binom{j}{k} a^{d-2k}=\frac{1}{2}
   \left(\left(-\sqrt{a^2-1}-a\right)^d+\left(\sqrt{a^2-1}-a\right)^d\right).
\end{equation*}}

 To prove that the absolute point is an $a'$-fold point of the curve $CH(n,d,a)$, it is sufficient to prove that any line passing through the absolute point intersects the curve at $2(n+d)$ points, $a'$ of which coincide with that point.
If we switch to the homogeneous coordinates ($x=\frac{x_1}{x_0}$, $y=\frac{x_2}{x_0}$), the above equations turn into:

\smallskip
\noindent -- if $n$ is an even number ($d$ must be odd)
{\small
\begin{align}
\label{implicit_n_par_hom}
& \hspace*{-0.5cm}  (x_1^2+x_2^2)^{n+1}\Big( \sum_{j=0}^{\lfloor \frac{d}{2}\rfloor}\sum_{k=0}^{j}\sum_{l=0}^{\lfloor \frac{d-2k}{2}\rfloor} (-1)^{d-k-1}\binom{d}{2j}\binom{j}{k}\binom{d-2k}{2l+1} a^{d-2k-2l-1}(x_1^2+x_2^2)^{l}x_0^{d-2l-1}\Big)^2\nonumber\\
=&\Big( \sum_{i=0}^{\lfloor \frac{n}{2}\rfloor} (-1)^i\binom{n}{2i} x_1^{n-2i}x_2^{2i}x_0^{d}\\
&-(x_1^2+x_2^2)^{\frac{n}{2}}\sum_{j=0}^{\lfloor \frac{d}{2}\rfloor}\sum_{k=0}^{j}\sum_{l=0}^{\lfloor \frac{d-2k}{2}\rfloor} (-1)^{d-k}\binom{d}{2j}\binom{j}{k}\binom{d-2k}{2l} a^{d-2k-2l}(x_1^2+x_2^2)^{l}x_0^{d-2l}\Big) ^2\nonumber
\end{align}}
-- if $n$ is an odd number
{\small
\begin{align}
\label{implicit_n_nep_hom}
& \hspace*{-0.5cm}  (x_1^2+x_2^2)^{n}\Big( \sum_{j=0}^{\lfloor \frac{d}{2}\rfloor}\sum_{k=0}^{j}\sum_{l=0}^{\lfloor \frac{d-2k}{2}\rfloor} (-1)^{d-k}\binom{d}{2j}\binom{j}{k}\binom{d-2k}{2l} a^{d-2k-2l}(x_1^2+x_2^2)^{l}x_0^{d-2l}\Big)^2\nonumber\\
=&\Big( \sum_{i=0}^{\lfloor \frac{n}{2}\rfloor} (-1)^i\binom{n}{2i} x_1^{n-2i}x_2^{2i}x_0^{d}\\
& - (x_1^2+x_2^2)^{\frac{n+1}{2}}\sum_{j=0}^{\lfloor \frac{d}{2}\rfloor}\sum_{k=0}^{j}\sum_{l=0}^{\mathrm D} (-1)^{d-k-1}\binom{d}{2j}\binom{j}{k}\binom{d-2k}{2l+1} a^{d-2k-2l-1}(x_1^2+x_2^2)^{l}x_0^{d-2l-1}\Big)^2\nonumber
\end{align}
}

Every straight  line $t$ passing through the absolute point $F_1(0,1,i)$ has the equation of the form
\begin{equation}\label{t}
x_2=ix_1+mx_0.
\end{equation}
It intersects the curve $CH(n,d,a)$ at the points whose coordinates  satisfy the equation obtained when $x_2$ from (\ref{t}) is substituted in (\ref{implicit_n_par_hom}) and (\ref{implicit_n_nep_hom}). Some short calculations deliver the following conclusions:
if $n>d$, for every line $t$, $x_0=0$ is the solution of the obtained equations with multiplicity $a'=2d$, and $F_1(0,1,i)$ is the intersection of $t$ and $CH(n,d,a)$ with the intersection multiplicity $a'=2d$.
If $d>n$, for every line $t$, $x_0=0$ is the solution  of the obtained equations with multiplicity $a'=n+d$, and $F_1(0,1,i)$ is the intersection of $t$ and $CH(n,d,a)$ with the intersection multiplicity $a'=n+d$.
To summarize: the absolute points are $2d$-fold points of $CH(n,d,a)$ when $d<n$, and they are $(n+d)$-fold points when $n<d$.
\hfill{$\square$}

\section{Generalized Rose Surfaces}

 In  paper \cite{circular} the authors considered a {\it congruence of circles} $\mathcal C(p)$ that consists of circles in Euclidean space $\mathbb E^3$  passing through two given points $P_1$, $P_2$. The points $P_1$, $P_2$ lie on the axis $z$ and are given by the coordinates $(0,0,\pm p)$, where $p=\sqrt{q}$, $q\in\mathbb R$. If $q$ is greater, equal to or less than zero, i.e. the points $P_1$, $P_2$ are real and different, coinciding (the axis $z$ is the tangent line of all circles of the congruence) or imaginary, the congruence $\mathcal C(p)$ is called an {\it elliptic},   {\it parabolic} or  {\it hyperbolic} congruence, respectively.

For every point $A$ $(A \notin z)$,  there exists a unique circle $c^A(p)\in\mathcal C(p)$ passing through the points $A$, $P_1$ and $P_2$.   For every point at infinity $A^\infty$ $(A^\infty\notin z)$, the circle $c^{A^\infty}(p)$ splits into the axis $z$ and the line at infinity in the plane $\zeta$ through the axis $z$ and the point $A^\infty$.

A point is defined as a  singular point of a congruence if infinitely many curves pass through it.   The singular points of $\mathcal C(p)$ are the points on the axis $z$ and the absolute points of $\mathbb E^3$ \cite{circular}.

For a given congruence $\mathcal C(p)$ and a given curve $\alpha$,  a {\it circular surface} $\mathcal {CS}(\alpha,p)$ is defined as the system of circles from $\mathcal C(p)$  that intersect $\alpha$.

\begin{center}
\includegraphics[scale=0.7]{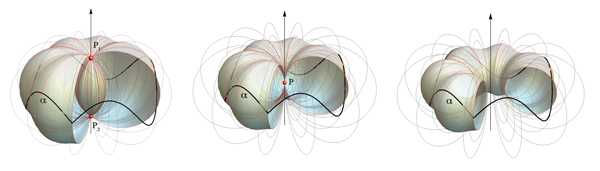}

\vspace{-0.5cm}

\begin{multicols}{3} a 

b 

c 
\end{multicols}
Figure 2: \parbox[t]{120mm}{Some examples of $\mathcal {CS}(\alpha,p)$ when  $\mathcal C(p)$ is an elliptic, parabolic and hyperbolic congruence are shown in figures a, b and c, respectively.}
\end{center}

If  $\alpha\colon I\rightarrow\mathbb R^3$, $I\subseteq \mathbb R$, is given by
$\alpha(t)=(\alpha_1(t), \alpha_2(t), \alpha_3(t))$, $\alpha_i \in C^1(I)$, and if its  normal projection  on the plane $z=0$ is denoted by $\alpha_{xy}(t)=(\alpha_1(t), \alpha_2(t), 0)$, then
 $\mathcal {CS}(\alpha,p)$ is given by the   parametric equations:
\begin{align}
\label{par.surface1}
x(t,\theta)&= \frac{\alpha_1(t)}{2\left \| \alpha_{xy}(t) \right \|^2}\left(\sqrt{4p^2\left \| \alpha_{xy}(t) \right \|^2 +(\left \| \alpha(t) \right \| ^2-p^2)^2}\cos\theta+\left \| \alpha(t) \right \| ^2-p^2\right)\nonumber\\
y(t,\theta)&= \frac{\alpha_2(t)}{2\left \| \alpha_{xy}(t) \right \|^2}\left(\sqrt{4p^2\left \| \alpha_{xy}(t) \right \|^2 +(\left \| \alpha(t) \right \| ^2-p^2)^2}\cos\theta+\left \| \alpha(t) \right \| ^2-p^2\right)\nonumber\\
z(t,\theta)&= \sqrt{p^2+\frac{(\left \| \alpha(t) \right \| -p^2)^2}{4\left \| \alpha_{xy}(t) \right \| ^2}}\sin\theta, \quad\quad (t, \theta) \in I \times [0,2\pi).
\end{align}

In \cite{circular} the authors showed that if $\alpha$ is an $m^{th}$ order algebraic curve that cuts the axis $z$ at $z'$ points, the absolute conic at $a'$ pairs of the absolute points and with the points $P_1$ and $P_2$ as $p'_1$-fold and $p'_2$-fold points, respectively, then, the following statements hold:

\linesp{0.9}
\begin{enumerate}
\item $\mathcal{CS}(\alpha,p)$ is an algebraic surface of the order $3m-(z'+2a'+2p'_1+2p'_2)$.
\item  The  absolute conic is an  $m-(z'+p'_1+p'_2)$-fold curve of $\mathcal{CS}(\alpha,p)$.
\item The axis $z$ is an $(m-2a'+z')$--fold line of $\mathcal{CS}(\alpha,p)$.
\item  The points  $P_1$, $P_2$ are  $2m-(2a'+p_1'+p_2')$--fold points of \,$\mathcal{CS}(\alpha,p)$.
\end{enumerate}

The singular lines and points of $\mathcal {CS}(\alpha,p)$ are the following:
\begin{itemize}
\item[5.] The line $z$ and the absolute conic are in the general case the singular lines of $\mathcal {CS}(\alpha,p)$. In particular, the points $P_1, P_2$   are of  higher multiplicity than the other points on the line $z$.
\item[6.]
 If different points of the curve $\alpha$ determine the same circle $c\in\mathcal C(p)$, this circle is a singular circle of  $\mathcal {CS}(\alpha,p)$ and its multiplicity equals the number of  intersection points of $\alpha$ and $c$. Through every singular point of $\alpha$, the singular circle of $\mathcal {CS}(\alpha,p)$ passes.
\item[7.] The intersection points  of $\alpha$ and $c(0)$, where  $c(0)$ is the circle given by the equations $x^2+y^2=-p^2$ $(p^2<0)$ and  $z=0$, are the singular points of $\mathcal {CS}(\alpha,p)$.
\end{itemize}

In this paper we study the properties of a special class of circular surfaces given by the following definition:

\begin{defi}
A generalized rose surface is a circular surface $\mathcal {CS}(\alpha,p)$ when the curve $\alpha$ is a generalized rose  $CH(n,d,a)$ lying in a plane that does not contain the axis $z$.
\end{defi}

Some generalized rose surfaces  $\mathcal {CS}(\alpha,p)$ for $\alpha =CH(n,d,0)$ and $p\in\mathbb R$ are simply called rose surfaces and were studied in  detail in \cite{roses}. Their parametric and implicit equations were derived, their singularities were investigated and they were visualized by the program {\it Mathematica}. This paper presents the extension of  \cite{roses} since we are also dealing with $\mathcal {CS}(\alpha,p)$ for $\alpha =CH(n,d,a)$ when $a \neq 0$ and $p\in\mathbb C$.

\bigskip

Let $CH(n,d,a)$ lie in any plane of $\mathbb E^3$ that does not contain the axis $z$, let its $2n$-fold ($n$-fold) point be denoted by $O$, and let its another real $j$-fold points be denoted by $D^j$.

 According to the properties 1 -- 4 of $\mathcal {CS}(\alpha,p)$  and the properties of $CH(n,d,a)$ given in Table 1, a  classification of generalized rose surfaces according to their order and the multiplicity of the axes $z$, the absolute points and the points $P_1$, $P_2$ can be made. Since this classification depends on the position of $CH(n,d,a)$ with respect to the singular points of the congruence $\mathcal {C}(p)$, the following cases should be taken into consideration:
\begin{itemize}
 \item[] {\bf Type 1}: $O=P_1$ or $O=P_2$;
\item[] {\bf Type 2}:  $O \in z$, $O \neq P_i$, $i=1, 2$;
\item[] {\bf Type 3}:  $D^j=P_1$ or $D^j=P_2$;
\item[] {\bf Type 4}:  $D^j \in z$, $D^j \neq P_i$, $i=1, 2$;
\item[] {\bf Type 5}:  $\forall Z \in z$, $Z \notin CH(n,d,a)$.
\end{itemize}
The algebraic properties of these surfaces are given in the following theorem:

\begin{tm}\label{generalized}

For generalized rose surfaces  $\mathcal {CS}(CH(n,d,a),p)$, the following table is valid:

\begin{center}
{\rm 

{\small

\begin{tabular}{|c|l|c| c| c| c |c |c |}
\hline
\multicolumn{3}{|c|}{ } &    &   &  & \\

\multicolumn{3}{|c|}{ } &     order of     &    multiplicity of     &    multiplicity of  &    multiplicity of \\
\multicolumn{3}{|c|}{  Type of $\mathcal{CS}(CH,p)$ } &    $\mathcal{CS}(CH,p)$    &   the absolute conic    &     the axis $z$  &     the points $P_i$\\
\multicolumn{3}{|c|}{ } &    &   &  & \\

\hline\hline
&$a=0$  &  $d<n$  &    $n+d$   &    $d$     &    $n-d$    &    $n$ \\

1A&$n\cdot d$   odd   &   $d>n$   &   $2d$   &   $d$    &   $0$   &    $d$ \\
\hline
&other  &   $d<n$   &   $2(n+d)$  & $2d$    &   $2(n-d)$   &  $2n$ \\

1B&cases   &   $d>n$   &    $4d$   &   $2d$    &   $0$    &   $2d$ \\
\hline\hline
&$a=0$   &   $d<n$   &    $2n+d$    &  $d$     &    $2n-d$   &   $2n$ \\

  2A &$n\cdot d$  odd   &   $d>n$  &  $n+2d$  &  $d$   &   $n$  &   $n+d$ \\
\hline
&other &    $d<n$   &  $2(2n+d)$  &   $2d$  &    $2(2n-d)$   &   $4n$ \\

  2B &cases  &   $d>n$   &    $2(n+2d)$   &   $2d$    &   $2n$   &   $2(n+d)$ \\
\hline\hline
&$a=0$  &  $d<n$   &   $3n+d-2j$   &   $n+d-j$     &   $n-d$  &   $2n-j$ \\

  3A &$n\cdot d$   odd  &  $d>n$   &   $2(n+d)-2j$  &  $n+d-j$   &  $0$    &    $n+d-j$ \\
\hline
&other  &    $d<n$  &  $2(3n+d)-2j$   &  $2(n+d)-j$   &    $2(n-d)$   &  $4n-j$ \\

  3B &cases   &    $d>n$  &    $4(n+d)-2j$  &  $2(n+d)-j$    &   $0$   & $2(n+d)-j$ \\
\hline\hline

&$a=0$  &  $d<n$  &    $3n+d-j$   &   $n+d-j$     &    $n-d+j$   &    $2n$ \\

4A &$n\cdot d$   odd   &   $d>n$  &  $2(n+d)-j$   &   $n+d-j$    &   $j$    &   $n+d$ \\
\hline
&other   &    $d<n$    &   $2(3n+d)-j$   &   $2(n+d)-j$    &    $2(n-d)+j$   &  $4n$ \\

 4B &cases   &   $d>n$   &   $4(n+d)-j$   &  $2(n+d)-j$    &   $j$    &   $2(n+d)$ \\
\hline\hline
&$a=0$  &  $d<n$   &    $3n+d$   &    $n+d$   &    $n-d$    &    $2n$ \\

 5A &$n\cdot d$   odd   &  $d>n$  &  $2(n+d)$   &   $n+d$   &   $0$   &    $n+d$ \\
\hline
&other   &   $d<n$   &   $2(3n+d)$   &  $2(n+d)$    &   $2(n-d)$   &  $4n$ \\

 5B &cases   &  $d>n$  &    $4(n+d)$   &   $2(n+d)$    &   $0$   &  $2(n+d)$ \\
\hline

\end{tabular}
}
\smallskip

Table 2: Properties of $\mathcal {CS}(\alpha,p)$

}
\end{center}

\end{tm}

\noindent {\bf Proof.} The statements regarding generalized rose surfaces are the immediate consequence of  the  algebraic properties 1 - 4 of the  circular surfaces $\mathcal {CS}(\alpha,p)$ when $\alpha$ is  the cyclic harmonic curve $CH(n,d,a)$ with the properties given in Table 1. More precisely, the values  $m$ and $a'$ are given in the second and fourth columns of Table 1, while  the values $p_1', p_2'$ and $z'$ depend on the position of $CH(n,d,a)$ with respect to the line $z$  as follows:\\
Type 1: Let us assume that $O=P_1$. Then $z'=p_2'=0$, while $p'_1=n$ when $a=0$ and $n \cdot d$ is odd, and $p'_1=2n$ in other cases.\\
Type 2:  $p_1'=p_2'=0$, $z'=n$  when $a=0$ and $n \cdot d$ is odd, and $z'=2n$ in other cases.\\
Type 3: If we assume that $D^j=P_1$, then $p_1'=j$, $p_2'=z'=0$.\\
Type 4:  $p_1'=p_2'=0$, $z'=j$.\\
 Type 5:  $p_1'=p_2'=z'=0$.

By substituting $m, a', p_1', p_2'$ and $z'$ into expression given in statements 1-4,  the algebraic properties of  $\mathcal {CS}(CH(n,d,a),p)$ are obtained. 
\hfill{$\square$}

\section{Vizualization of Generalized Rose Surfaces}

In this section we give some examples of  generalized rose surfaces. For a particular surface,   parametric and implicit equations can be obtained on the basis of the equation (\ref{par.surface1}). We check the properties given in Theorem \ref{generalized} and visualize the shape of the surface. For computing and plotting we use  the program {\it Mathematica}.

\bigskip

The circular surfaces $\mathcal {CS}(CH,0)$ for $CH=CH(7,3,a)$, where $a=0, 0.25, 1$ and $2.5$, are shown in Figures  3 a, b, c and d, respectively. The curve $CH$ lies in the plane $z=0$ and the point $O$ coincides with the directing point $P_1=P_2$ of the parabolic  congruences  $\mathcal C(0)$. Therefore, $CH$ in the case a is an algebraic curve of the order 10 having the absolute conic as a triple curve, axis $z$ as a quadruple line and point $O$ as a 14-fold point. The other three surfaces are of the order 20 and  the absolute conic is a 6-fold curve, the axis $z$ is an 8-fold line and the point $O$ is a 14-fold point of the surfaces.

\vspace{-0.5cm}

\begin{center}
\includegraphics[scale=0.6]{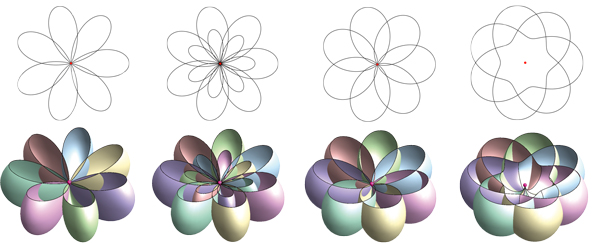}

\vspace{-0.5cm}

\begin{multicols}{4} a 

b 

c

d
\end{multicols}

Figure 3: Four surfaces $\mathcal {CS}(CH,0)$ of type 1

\end{center}

In Figure 4 three generalized rose surfaces  $\mathcal {CS}(CH,1)$  of type 1B are presented.  $CH$ in Figures  4 a, b and c are curtates $CH(3,1,1.25)$, $CH(2,3,1.25)$ and $CH(7,3,1.25)$ in the plane $z=-1$ with  the $2n$-fold isolated point $O$ coinciding with the directing point $P_1$ of the elliptic  congruence  $\mathcal C(1)$, respectively.

\begin{center}
\includegraphics[scale=0.65]{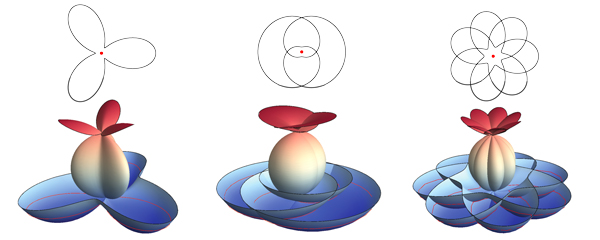}

\vspace{-0.5cm}

\begin{multicols}{3} a 

b 

c 
\end{multicols}

Figure 4: Three surfaces $\mathcal {CS}(CH,1)$ of type 1B

\end{center}

If a circular surface is defined by a hyperbolic congruence and a curtate with the $2n$-fold point on the axis $z$,  the surface is of type 2B. The generalized  rose surfaces  $\mathcal {CS}(CH,i)$ for $CH(9,2,2)$ in the plane $z=0, 0.5, 1$ are shown in Figures 5 a, b, c, respectively.  All three  surfaces  $\mathcal {CS}(CH,i)$  are of the order 40. The absolute conic is a quadruple  curve and the axis $z$ is a 32-fold line of the surfaces.

\begin{center}
\includegraphics[scale=0.7]{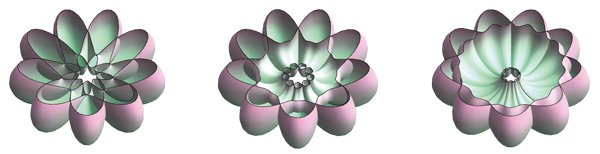}

\vspace{-0.5cm}

\begin{multicols}{3} a 

b 

c 
\end{multicols}

Figure 5: Three surfaces $\mathcal {CS}(CH,i)$ of type 2B

\end{center}

The surfaces   $\mathcal {CS}(CH,i)$  for $CH(7,1,2)$ in the plane $z=0.75, 0$ are shown in Figures 6 a and b,  respectively. The surface in Figure 6 c is defined by $CH(7,1,1.5)$ in the plane $z=0$. All three surfaces are of type 2B and of the order 30, with the absolute conic as a double curve, axis $z$ as a 26-fold line and the points $P_1, P_2$ as the 28-fold points. The surfaces in Figures 6 b and c possess the  singular points in the intersection points  of $CH$ and $c(0)$, where  $c(0)$ is the circle given by the equations $x^2+y^2=1$, $z=0$.

\begin{center}
\includegraphics[scale=0.7]{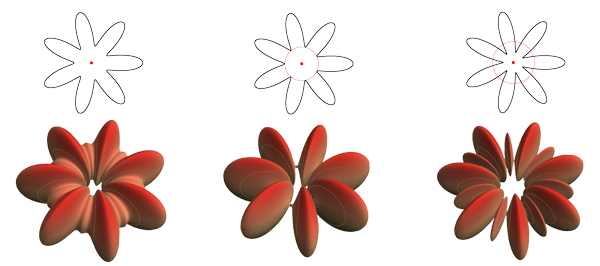}

\vspace{-0.5cm}

\begin{multicols}{3} a 

b 

c 
\end{multicols}

Figure 6: Three surfaces $\mathcal {CS}(CH,i)$ of type 2B

\end{center}

The surface  $\mathcal {CS}(CH,0)$  shown in Figure 7 a is defined
by the parabolic congruence $\mathcal C(0)$ and the rose $CH(3,1,0)$
intersecting the axis $z$ at its regular point $D^1=P_1=P_2$. The
obtained surface is of type 3A and the order 8. The absolute conic is a
triple conic, the axis $z$ is a double line and $D^1$  is an
eightfold point of the surface. The surfaces $\mathcal {CS}(CH,0)$
in Figures 7 b and c are of type 3B. The $CH(3,2,0.5)$ in Figure 7 b
has a triple point at the directing point of the congruence
$\mathcal C(0)$, while  $CH(3,2,0)$ in  Figure 7 c has a regular
point at that point.

\begin{center}
\includegraphics[scale=0.7]{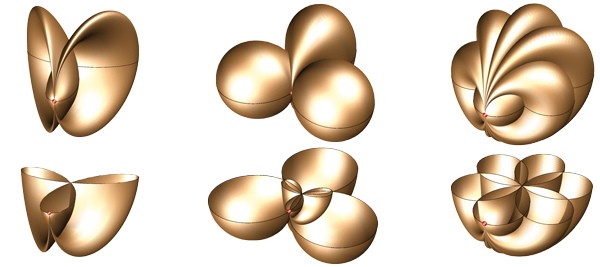}

\vspace{-0.5cm}

\begin{multicols}{3} a 

b 

c 
\end{multicols}

Figure 7: Three surfaces $\mathcal {CS}(CH,0)$ of type 3

\end{center}

All surfaces of type 4 pass through the real lines at infinity. Three such surfaces are presented in Figure 8.  $\mathcal {CS}(CH,p)$ in  case a is defined by  $CH(3,1,0)$, and $p=1$, $\mathcal {CS}(CH,p)$ in  case b is defined by  $CH(3,1,0)$ and $p=i$, and $\mathcal {CS}(CH,p)$ in case c is defined by  $CH(5,1,0)$ and $p=i$. All three roses lie in the plane $z=0$ and intersect the axis $z$ at the regular point $D^1$. The first two surfaces are of the order 9 with the absolute conic as a triple curve, axis $z$ as a triple line and the directing points as the sixfold points. The third surface is of the order 15 with the absolute conic as a 5-fold curve, axis $z$ as a 5-fold line and the directing points as the 10-fold points.

\begin{center}
\includegraphics[scale=0.7]{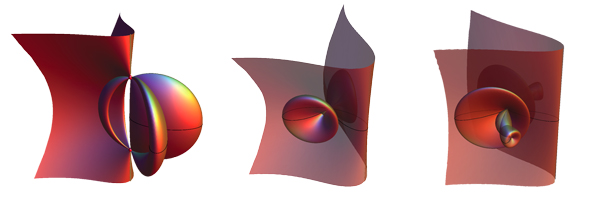}

\vspace{-0.5cm}

\begin{multicols}{3} a 

b 

c 
\end{multicols}

Figure 8: Three surfaces $\mathcal {CS}(CH,p)$ of type 4A

\end{center}

Three surfaces of type 5A defined by $CH(3,1,0)$ and $p=1,0,i$ are shown in Figure 9 a, b, c, respectively. $CH(3,1,0)$ lies in the plane $z=0$ and has a triple point at the point $(1,0,0)$. All three surfaces are of the order 10, have the absolute conic as a quadruple curve, axis $z$ as a double line and the directing points as the sixfold points.

\begin{center}
\includegraphics[scale=0.7]{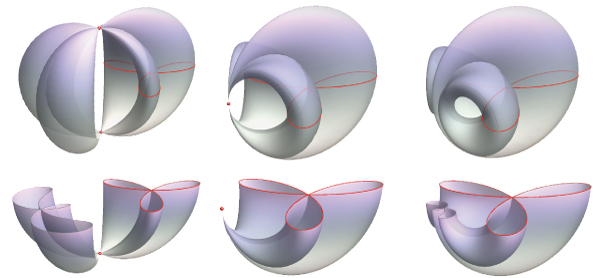}
\vspace{-0.5cm}

\begin{multicols}{3} a 

b 

c 
\end{multicols}

Figure 9: Three surfaces $\mathcal {CS}(CH,p)$ of type 5A

\end{center}


\end{document}